\magnification=\magstep1
\input amstex
\input epsf
\documentstyle{amsppt}
\voffset=-3pc
\loadmsbm
\parskip=6pt
\NoBlackBoxes
\def\ds{\displaystyle}
\def\og{\overline g}
\def\ox{\overline x}
\def\bZ{\Bbb Z}
\def\bC{\Bbb C}

\def\cZ{\Cal Z}
\def\cS{\Cal S}
\def\Res{\operatorname{Res}}

\def\var{\Cal E}
\def\vare{\varepsilon}
\def\uG{\bold G}
\def\uM{\bold M}
\def\uH{\bold H}
\def\uB{\bold B}
\def\uA{\bold A}
\def\uP{\bold P}

\def\wt{\widetilde}

\def\tP{\widetilde P}
\def\tM{\widetilde M}
\def\tN{\widetilde N}
\def\tG{\widetilde G}

\def\tA{\widetilde A}
\def\a{\alpha}
\def\O{\Omega}
\def\diag{\operatorname{diag}}
\def\Ind{\operatorname{Ind}}
\def\wt{\widetilde}
\topmatter
\title Reducibility for $SU_n$ and generic elliptic representations\endtitle
\author David Goldberg\endauthor
\endtopmatter

\subheading{Introduction} The problem of classifying the tempered spectrum of a connected reductive quasi-split group,
defined over a local field, $F,$ of characteristic zero, consists of three parts.  The first is to classify the discrete
series representations of any Levi subgroup.  The second step is to understand the rank one Plancherel measures, which is
equivalent to understanding the reducibility of those representations parabolically induced from a discrete series of maximal
Levi subgroups.  The third step is to understand the structure of representations parabolically induced from discrete
series representations of an arbitrary parabolic subgroup, using the second step and the combinatorial theory of the
Knapp-Stein
$R$--group.  We address this third step here for the case where the group in question is the quasi-split special unitary
group. 

This work builds on earlier results on $R$-groups due to several authors.  For quasi-split unitary groups, $U_n(F),$ the
$R$-groups attached to the principal series were computed by Keys \cite{\bf K1, K2} and the theory of restriction
yielded
$R$-groups for
$SU_n(F).$  Keys gave a description of the $R$--groups in terms of the Langlands parameterization, which in that case is well
understood \cite{\bf L}.  Recent work of Ban and Zhang \cite{\bf B-Z} has shown that the construction of the
$R$--group  through $L$--group considerations, as described in \cite{\bf A1}, will be valid for all quasi-split
connected reductive groups, once the local parameterization conjecture is established.  For
$U_n(F)$ and an arbitrary parabolic subgroup, the $R$--groups were computed in \cite{{\bf G2}} and these were
shown to be elementary two groups.  The computation of the $R$--group via the $L$--group was carried out by D. Prasad
\cite{\bf P},
and one can see that this computation agrees with \cite{\bf G2}, if you assume that the Langlands
$L$--functions are the same as the Artin $L$--functions attached to the parameter.  This has been established in some cases
by Henniart \cite{\bf H}.

Following the methods of \cite{\bf K1,K2}, and using the theory of \cite{\bf G-K,T},
we work via restriction.  As in \cite{\bf T} some of our results can
be proved in a more general setting.  Namely, we consider the case where $G\subset \widetilde G$ have the same derived group.
We are able to give a rough description of the $R$-group in this setting.  In the case of special unitary groups we make
this explicit due to a stronger understanding of the restriction of irreducible smooth representations from Levi
subgroups
$\widetilde M$ of $\widetilde G,$ to the subgroup $M=\widetilde M\cap G.$ These results on restriction are given in the
latter part of Section 2. 

To be more precise we consider a $p$--adic field $F$ of characteristic zero, and residual characteristic $q_F,$
and fix a quadratic extension
$E/F.$   Let
$\wt {\bold G}$ be a quasi-split reductive group defined over F, with $\wt{\bold G}(F)=\wt G.$  We assume that $\bold
G\subset
\wt{\bold G}$ is a reductive subgroup with the property that $[\wt G,\wt G]=[G,G],$ where these represent the derived
groups. Then, for any Levi subgroup $\bold{\wt M}$ of $\bold{\wt G},$ we have $MZ(\wt M)$ is of finite index in $\wt
M.$  Here $M=\bold M(F),$ with $\bold M=\wt{\bold M}\cap \bold G,$ and $Z(\wt M)$ is the center of $\wt M=\bold{\wt
M}(F).$  Thus, the theory of Section 2 of \cite{{\bf T}} applies.  
In particular, if $\sigma$ is a discrete series
representation of
$M,$ then there is a discrete series of $\wt M$ with $\sigma$ a component of $\pi|_{M}.$  Furthermore, any irreducible
representation $\pi'$ of
$\wt M$ for which
$\operatorname{Hom}_M(\sigma,\pi')\neq0$ is of the form $\pi'\simeq\pi\chi,$ for a character $\chi$ of $\wt M$ trivial on $M.$

The $R$-group, $R(\sigma),$ is a subgroup of the stabilizer, $W(\sigma),$ of $\sigma$ in the Weyl group $W,$  and by the
considerations above,
$$W(\sigma)\subset\{w|\pi^w\simeq\pi\chi,\text{ for some }\chi\}.$$  
The theory of the $R$--group dictates that the complement of $R(\sigma)$  in $W(\sigma)$ is the subgroup, $W'(\sigma),$
generated by the reflections
$w_{\alpha}$ for which the rank one Plancherel measure $\mu_\alpha(\sigma)=0.$  We show that $\mu_{\alpha}(\sigma)=0$ if and
only if $\mu_\alpha(\pi)=0,$ and thus $W'(\sigma)=W'(\pi)$ (cf. Lemma 2.3). 
 Let 
$$\widehat {W(\sigma)}=\{\chi|\pi^w\simeq\pi\chi\text{ for some }w\in W(\sigma)\}.$$
We show that, for any $\chi\in\widehat{W(\sigma)},$ there is some $w_{\chi}\in R(\sigma)$
with $\pi^{w_\chi}\simeq\pi\chi,$ and this element is unique up to multiplication by elements of $R(\pi)\cap W(\sigma).$
Thus, in this general setting we always have $R(\sigma)/\left(R\left(\pi\right)\cap
W(\sigma)\right)\simeq \widehat{W(\sigma)}/X(\pi),$ where $X(\pi)=\{\chi|\pi\chi\simeq\pi\}$ (cf Proposition 3.2).

We then specialize to the case where $\tilde{\bold G}=U_n$ and $\bold G=SU_n.$
There we show $R(\pi)\triangleleft R(\sigma).$  Furthermore, we can split the sequence $1\rightarrow R(\pi)\rightarrow
R(\sigma)\rightarrow \widehat{W(\sigma)}/X(\pi)\rightarrow1.$  Thus, there is a subgroup $\Gamma_{\sigma}$ of $R(\sigma)$
which is isomorphic to $\widehat{W(\sigma)}/X(\pi)$ for which $R(\sigma)=\Gamma_\sigma\ltimes R(\pi).$

The structure of $\operatorname{Ind}_P^G(\sigma)$ is determined by the representation theory of a certain extension of
$R(\sigma).$  More precisely, for $w\in R(\sigma),$ we choose an intertwining operator satisfying $T_w\sigma^w=\sigma T_w.$
Then there is a $2$--cocycle $\gamma: R(\sigma)\times R(\sigma)\rightarrow \Bbb C$ defined by
$T_{w_1w_2}=\gamma(w_1,w_2)T_{w_1}T_{w_2}.$  The intertwining $\Cal C(\sigma)$ of $\operatorname{Ind}_P^G(\sigma)$ is
then isomorphic to the twisted group algebra $\Bbb C[R(\sigma]_\gamma.$  For simplicity of this exposition, we assume
that the cocycle splits (which is the case whenever $\sigma$ is generic \cite{\bf K2}, and is always the case for
for the classical groups for which $R$--groups have been computed \cite{\bf Hr,G2, G4}).  We note a correction to the description of the elliptic
tempered spectrum of $U(n)$ in \cite{\bf G2} (cf remark 3.10).
 Then there is a
correspondence
$\rho\mapsto\pi(\rho)$ between irreducible representations of
$R(\sigma)$ and classes of irreducible components of $i_{G,M}(\sigma)$ \cite{\bf A2, K2}. The multiplicity of
$\pi(\rho)$ in
$i_{G,M}(\sigma)$ is equal to $\dim\rho.$  Further, the behavior of the character $\theta_\rho$ of $\rho$
 determines which components of $i_{G,M}(\sigma)$ are elliptic \cite{\bf A2}.   More precisely, if $R(\sigma)_{reg}$ is the set of
elements of $R(\sigma)$ for which the fixed points $\frak a_w$ of $w$ in $\frak a,$ the real lie algebra of
$\bold A,$ is as small as possible, namely $\frak a_G,$  then $\pi_\rho$ is elliptic if and only if $\theta_\rho$ is non-vanishing on
the regular set $R(\sigma)_{reg}.$ In our situation, we examine the elliptic spectrum in the case where
$\pi$ is generic. recall that the Weyl group is the semidirect product of a permutation group with an elementary two group (consisting of ``sign changes''). We show that 
$i_{G,M}(\sigma)$ has an elliptic component, if and only if there is an element $w=sc$  of $R(\sigma)$ whose permutation component $s$
is of maximal possible length, and if $c$ changes an odd number of signs.  Finally we give an example of a phenomenon which we had not
noted before.  Namely, a case where $i_{G,M}(\sigma)$ has some  elliptic components, but not all of the components are elliptic.  In
fact, induced representations which have this property were exhibited in \cite{\bf K2}, but, as this predated Arthur's description
of the elliptic spectrum, \cite{\bf A2}, it was not noted there.  We give a specific new example, and indicate how this can be generalized.

I wish to thank J. Arthur, D. Ban, G. Henniart, R. Kottwitz, F. Shahidi, and M. Tadic for helpful conversations and communications
regarding this work.  Finally, this work grew out of conversations with Alan Roche regarding the connection  between Arthur's
$R$--group construction and the Knapp-Stein $R$--group.  I wish to  especially thank Alan for his insights, comments, and
collegial support.

\subheading{\S1 Notation and preliminaries}

Let $F$ be a nonarchimedean local field of characteristic zero, and residual characteristic $q_F.$ Fix a quadratic
extension
$E/F$. Let
$\gamma$ be the non--trivial Galois automorphism of $E/F$, which we also denote by $x\mapsto\ox.$ Fix
an element $\beta\in E$ with
$\gamma(\beta)=-\beta$. For $n\in\bZ^+$, let $u_n=\pmatrix &&&&&.\\ &&&&.\\&&&.\\ &&-1\\ &1\\ -1\endpmatrix$, and fix a hermitian form
$h_n\in M_n(E)$, by
$$ h_n=\cases u_n&\text{if $n$ is odd},\\
\beta u_n&\text{if $n$ is even.}\endcases
$$ For $g=(g_{ij})\in\Res_{E/F}GL_n$, we let $\og=(\og_{ij})$ and set $\vare(g)=u_n\,^t\bar g^{-1} u_n^{-1}$. Denote by
$\widetilde{\uG}=\widetilde{\uG}(n)=U_n$, the quasi-split unitary group defined with respect to $E/F$ and $h_n$. Thus
$$
\widetilde{\uG}=\{g\in\Res_{E/F}GL_n|\, g h_n \,^t\og=h_n\},
$$

We let $\uG=\uG(n)=SU_n=\bold\tG\cap\Res_{E/F}SL_n$. If $\widetilde{\uH}\subset\widetilde{\uG}$, then we let $\uH=\widetilde{\uH}\cap\uG$.
We denote the $F$--points by $\tG=\widetilde{\uG}(F)$ and $G=\uG(F)$, and similarly for other groups.
Let $\widetilde{\bold T}$ be the maximal torus in $\widetilde{\uG}$ of diagonal elements, and let $\widetilde{\uA_0}$
be the maximal split subtorus of $\widetilde{\bold T}.$  Denote by $\Phi(\widetilde{\uG},\widetilde{\uA}_0)$ the roots
of $\widetilde{\uA}_0$ in $\widetilde{\uG}.$  We fix $\widetilde{\uB}=
\widetilde{\bold T}\widetilde{\bold U}$ to be the
Borel subgroup of upper triangular elements of $\widetilde{\uG}.$  Let $\widetilde\Delta$ be the simple roots with
respect to this fixed choice of Borel subgroup.  If $\theta\subset\widetilde\Delta,$ then we denote by
$\bold{\tP}_\theta=\bold{\tM}_\theta\bold{\tN}_\theta$ the associated standard parabolic subgroup. Since we will
be working with a fixed $\theta,$ we drop the subscript and simply write $\bold{\tP}=\bold{\tM\tN}.$ 
Let $\wt{\uA_\uM}$ be the split component of $\uM.$
 We write
$H_{\tilde P}$ (respectively, $H_P$) for the homomorphism from $\tilde M$ (respectively M) to $\frak a_{\tilde M}$
(respectively $\frak a_M)$ given in
\cite{\bf HC}, where
$\frak a_{\tilde M}$ and $\frak a_M$ are the real lie algebras of $A_{\tilde M}$ and $A_M,$ respectively.

Note that, for some choice of partition $\{n_1,n_2,\dots,n_r,m\}$ of $\ds{\left[\frac{n}{2}\right],}$
$$\wt{\uA_\uM}=\left\{\diag\{x_1I_{n_1},x_2I_{n_2},\dots,x_rI_{n_r},I_m,\bar x_r^{-1}I_{n_r},\dots,\bar x_1^{-1}
I_{n_1}\}|x_i\in\Res_{E/F} GL_1\right\},$$
and
$$
\gather
\wt{\uM}=\left\{\diag\left\{g_1,g_2,\dots,g_r,h,\vare(g_r),\dots,\vare(g_1)\right\}\,|g_i\in\Res_{E/F} GL_{n_i},h\in  U(m)\right\}\simeq\\
\Res_{E/F} GL_{n_1}\times\Res_{E/F}GL_{n_2}\times\dots\times \Res_{E/F} GL_{n_r}\times U(n).
\endgather$$

Let $\Phi(\bold\tG,\bold\tA)$ be the reduced roots of  $\widetilde{\bold A}$ in $\bold\tG.$
The
Weyl group, $W(\widetilde{\uM})=N_{\widetilde{\uG}}(\widetilde{\uA})/\widetilde{\uM}
\simeq\cS\rtimes\bZ_2^r,$ where $\cS\subset S_r$ is generated by the transpositions $(ij)$ for which $n_i=n_j$. The realization of
$s_{ij}\in W(\widetilde{\bold M})$ with $s_{ij}\mapsto (ij)$ under the above isomorphism is given by
$$
\gathered s_{ij}(g_1,\ldots g_i,\ldots,g_j,\ldots g_r,h)=\\ (g_1,\ldots g_j,\ldots g_i,\ldots g_r,h).
\endgathered
$$ The subgroup $\cZ\simeq\bZ_2^r$ of $W(\widetilde{\uM})$ is generated by ``block sign changes'' $C_i$ given by
$$ C_i(g_1,g_2,\ldots g_i,\ldots g_r,h)=(g_1,g_2,\ldots\vare(g_i),\ldots g_r,h).
$$ 
We use $w$ to represent both a class in $W(\widetilde M)$ and a representative of that class in
$N_G(\widetilde M)$. This should cause no confusion here.  If $\bold M=\widetilde{\bold M}\cap \bold G,$ and $\bold
A=\widetilde{\bold A}\cap \bold G,$
then $\Phi(\widetilde{\uG},\widetilde{\uA})=\Phi(\uG,\uA),$ and $W(\uM)=W(\widetilde{\uM}),$ and we identify the action
of $W(\uM)$ as the restriction of the action of $W(\widetilde{\uM})$ to $\uM.$

If $\pi$ is an irreducible admissible representation of $\tM$, then
$\pi\simeq\pi_1\otimes\ldots\otimes\pi_r\otimes\tau$, where each
$\pi_i$ is an irreducible admissible representation of $GL_{n_i}(E)$ and $\tau$ is one of $G(m)$.
Then $W=W(\widetilde{\uM})$ acts on $\pi$ by $\pi^w(m)=\pi(w^{-1}mw)$. 

We use Harish--Chandra's notation of $\var_c(\tG)$ to denote the (equivalence classes of) irreducible admissible representations of $\tG,$
and  $\var_t(\tG),\ \var_2(\tG),$ represent the tempered, and square integrable classes, respectively.  Similar
notation is used for
$\tM,\,M,$ and $G.$  We use the notation $i_{G,M}(\sigma)$ for the representation of $G$ unitarily induced from the representation
$\sigma$ of $M,$ of course extended trivially to the unipotent radical of a parabolic $P$ with Levi component $M.$

Note that the center $Z(\tG)$ of $\tG$ is isomorphic to $E^1=\{x\in E|x\overline x=1\}$. Since $G=[\tG, \tG]$, the
theory of Section 2 of \cite{\bf T} applies to restriction from $\tG$ to $G$. Similarly we can apply these results to
restriction from
$\tM$ to
$M$, which is the subject of the next section.

We use $X(\tG)$ and $X(\tM)$ to denote the characters of $\tG$ and $\tM$, respectively.  If $\chi\in X(\widetilde G),$ then
$\chi(g)=\chi'(\det g),$ for some character $\chi'$ of  $E^1.$  We will abuse notation and use $\chi$ to represent
both the character of $\wt G,$ and the character of $E^1$ through which it factors.   If $\pi\in\var_2(\wt M),$ and
$\chi\in X(\wt M)$ then we denote by $\pi\chi$ the representation $g\mapsto\pi(g)\chi(g).$  Then we let
$X(\pi)=\{\chi\in X(\wt M)\, |\pi\chi\simeq\pi\}.$

In order to describe the generic elliptic spectrum we need to make some observations about the Lie algebra $\frak g$
of $G.$  We have 
$$\frak g\simeq\left\{X\in\frak sl_{n}(\Bbb  C)| XJ+ J\,^ t\bar X=0\right\}.$$
If $\frak a_M$ is the  Lie algebra of $A_M,$ then a straightforward calculation shows
$$\frak a_M=\{\diag(a_1I_{n_1},a_2I_{n_2},\dots, a_r I_{n_r},0_m,-\bar a_{r} I_{n_r},\dots, -\bar a_1
I_{n_1})\,\big|\sum_i n_i(a_i-\bar a_i)=0\}.$$ We will sometimes denote an element of $\frak a_M$ by $Y_M(a_1,
\dots,a_r).$  For $w\in W(\bold G, \bold A_M),$ we let $\frak a_w=\{X\in\frak a_M|w\cdot
X=X\},$ where $w\cdot X=ad(w)X.$ Note that $\frak a_G=\{0\}.$

\subheading{\S2 Plancherel measures and restriction}

In this section we establish some results on compatibility of Plancherel measures with restriction. Much of this follows immediately from
earlier work to which we refer. Many of 
these results apply more generally to the situation where $G\subset\widetilde G$ and
$[G, G]=[\widetilde G,\widetilde G].$ For the moment we work in that context,
with $\tM$ be a Levi subgroup of $\tG$ (possibly
$\tG$ itself) and $M=\tM\cap G.$

\proclaim{Lemma 2.1}Let $\pi\in\var_2(\tM)$.

\roster
\item"(a)" There is an integer $m_0$ so that $\pi|_M=m_0\ds\bigoplus^k_{i=1}\sigma_i$, with $\sigma_i$ irreducible and inequivalent.
\cite{\bf G-K}.

\item"(b)" If Hom$_M(\pi,\pi')\not= 0$, then there is a character $\chi$ of $E^1$, so that $\pi'\simeq\pi\chi$. 
\cite{\bf T, Corollary
2.5}.

\item"(c)" 
Every $\sigma\in\var_2(M)$ is a component of $\pi|_M$ for some $\pi\in\var_2(\tM)$ \cite{\bf T, Prop.
2.2}.\qed
\endroster
\endproclaim

\proclaim{Corollary 2.2}Suppose $\sigma\in\var_2(M)$ and $\pi\in\var_2(\tM)$ with $\pi|_M$ containing $\sigma$. Suppose that, for some
$w\in W(\uM)$, $\sigma^w\simeq\sigma$. Then, there is a $\chi\in X(\tM)$ so that $\pi^w\simeq\pi\chi$.\qed
\endproclaim


For a reduced root $\alpha\in\Phi(\widetilde{\uP},\widetilde{\uA})$, we write $\mu_\alpha(\pi)$ for the rank one Plancherel measure
attached to $\alpha$ (see \cite{\bf H-C}). We know that $\mu_\alpha(\pi)=0$ if and only if the standard intertwining
operator
$\nu\mapsto A(\nu,\pi,w_\alpha)$ has a pole at $\nu=0$. Here $\nu\in\wt{\frak a}^*_{\bC}$, the complexified dual of
the Lie algebra of
$\widetilde{\uA}$, and $w_\alpha$ is the reflection associated to $\alpha$.  Similarly, for $\sigma\in\var_2(M)$ we have
the Plancherel measure
$\mu_\alpha(\sigma)$. This is given by the pole of the standard intertwining operator $\nu_0\mapsto
A(\nu_0,\sigma,\mu_\alpha),$ where $\nu_0\in\frak a^*,$ the complexified dual of the Lie
algebra of $\uA.$
Note
that, since $\tN=N$, the intertwining operators are given by the same formula. That is, for $\nu$ and $\nu_0$ in
the region of convergence,
$$
\align A(\nu,\pi,w_\alpha)\widetilde f(\widetilde g)&=\int_{\,^*\!N_\alpha}\widetilde f(w_\alpha^{-1} n \widetilde g) dn,\text{
and}\tag2.1\\ A(\nu_0,\sigma,w_\alpha) f(g)&=\int_{^*\!N_\alpha} f(w_\alpha^{-1} n g) dn,\tag2.2
\endalign
$$ Where $\widetilde f$ is in the space of $i_{\widetilde M_\alpha,\widetilde M}(\pi\otimes q_F^{<\nu, H_{\tilde
P}()>}),\ f$ is in the space of
$i_{M_\alpha,M}(\sigma\otimes q_F^{<\nu_0,H_P()>}),\
\widetilde g\in\tM_\a$, and $g\in M_\alpha$ (see \cite{\bf H-C} for definitions of $M_\alpha,\ ^*N_\alpha$).  The
operators are then defined for all $\nu$ by meromorphic continuation.

\proclaim{Lemma 2.3}Suppose $[\widetilde G,\widetilde G]=[G, G]$ and $ZG\backslash\widetilde G$ is finite
abelian. Let
$\widetilde M$ be a Levi subgroup of $\tG$ and set $M=\tM\cap G.$ Let $\pi\in\var_2(\tM)$ and suppose $\sigma\in\var_2(M)$ is a
component of
$\pi|_M$.  Then, for any $\alpha\in\Phi( \bold{ \tilde P}\bold{\tilde A})=\Phi(\bold P,\bold A),$ we have
$\mu_\alpha(\pi)=0$ if and only if $\mu_\alpha(\sigma)=0.$
\endproclaim

\demo{Proof}  We note that one can adapt the proof of a similar result in \cite{\bf Sh3}, but we choose a slightly different approach.
Let  $\pi|_M=m\bigoplus\limits^k_{i=1}\sigma_i$, and assume $\sigma=\sigma_1.$   If $\Pi=i_{\tilde G,\tilde M}(\pi),$
 then $\Pi|_{G}\simeq i_{G,M}(\pi|_{M})=m\oplus i_{G,M}(\sigma_i).$
Let $w=w_\alpha.$ 
Then the  intertwining operators for  $\tilde G$ satisfy 
$$A(w\nu,w\pi,w^{-1})A(\nu,\pi,w)=\mu_\alpha(\nu,\pi)^{-1}\cdot\gamma_{\alpha}(\tilde G/\tilde P)^2,\tag 2.3$$
where 
$$\gamma_{\alpha}(\tilde G/\tilde P)=\int_{^*N_\a}q_F^{-2<\rho_\alpha,H_{P_{\a}}(n)>}\,dn,$$ 
is the constant given in {\cite{\bf H-C}.
The meromorphic function $\mu_{\alpha}(\nu,\pi)$ is the Plancherel density and $\mu_\alpha(0,\pi)=\mu_\alpha(\pi).$
Restricting the relation (2.3) to $G,$ (and restricting to the $i_{G,M}(\sigma\otimes q_F^{<\nu|_{\frak
a},-->})$-isotypic subspace)  gives
$$A(w\nu|_{\frak a},w\sigma,w^{-1})A(\nu|_{\frak a},\sigma,w)=\mu_\alpha(\nu|_{\frak a},\pi)^{-1}\cdot\gamma_{\alpha}(
\tilde G/ \tilde P)^2.\tag 2.4$$ On the other hand, we have
$$A(w\nu|_{\frak a},w\sigma,w^{-1})A(\nu|_{\frak a},\sigma,w)=\mu_\alpha(\nu|_{\frak
a} ,\sigma)^{-1}\cdot\gamma_{\alpha}( G/ P)^2,\tag 2.5$$
and the result follows by letting $\nu$ go to zero.
\qed
\enddemo

Let $\Delta'(\pi)=\{\alpha\in\Phi(\wt{\bold P},\wt{\bold A})|\mu_\alpha(\pi)=0\}.$  Lemma 2.3  shows
$\Delta'(\pi)=\Delta'(\sigma),$ where $\Delta'(\sigma)$ is similarly defined.
We let $W'(\pi)=<w_\alpha|\alpha\in\Delta'(\pi)>.$  Then
 $W'(\pi)=W'(\sigma),$ and since $\mu_\alpha(\sigma)=0$ implies $\sigma^{w_\alpha}\simeq\sigma,$ 
$W'(\pi)\subset
W(\sigma)=\{w\in W(M)|\sigma^w\simeq\sigma\}.$

We now give several partial  converses to Corollary 2.2 in the case where $\widetilde{\uG}=U_n$ and $\uG=SU_n.$ Some of these will be crucial
in describing the $R$-group explicitly in Section 3, while others we include to show the extent to which we can establish the converse at this
time.  We have yet to determine whether the converse holds in general.

\proclaim{Lemma 2.4} Suppose  $\wt M$ is maximal,  $\pi\in\var_c(\wt M)$ is generic, and 
$\pi^w\simeq\pi.$ Then, for each
component $\sigma$ of $\pi|_M$ we have $\sigma^w\simeq\sigma.$
\endproclaim
\demo{Proof} 
Let $\eta=\eta_{E/F}$ be the quadratic character of $F^\times$ attached to the extension $E/F$ by local class field theory.
Fix any character $\chi_\eta$ of $E^\times $ so that $\chi_\eta|_{F^\times}=\eta.$ 
First suppose $m=0,$ which implies $n=2k$ is even. Note that 
$$M=\left\{\pmatrix g\\&\vare(g)\endpmatrix\,\Bigg|\det g\det\vare(g)=1\right\}.$$
Since $\det\vare(g)=\overline{\det g^{-1}},$ we have $M\simeq\{g\in GL_k(E)|\det g\in F^\times\}.$ Thus,
$\chi_\eta|_{M}=\eta.$  Since $\wt M\simeq GL_k(E)$ and  $\pi^w\simeq\pi^\vare,$ we have
$\pi^\vare\simeq\pi.$ By \cite{\bf G3} exactly one of $L(s,\pi,A)$ or $L(s,\pi\chi_\eta,A)$ has a pole at $s=0,$
where $A$ is the Asai representation of  $GL_k(\Bbb C).$ Thus, either $\mu(\pi)=0$ or $\mu(\pi\chi_\eta)=0,$
where $\mu$ is the Plancherel measure \cite{\bf Sh2}.
Let $\sigma$ be an irreducible component of $\pi|M.$ 
Then, by Lemma 2.3 either $\mu(\sigma)=0$ or $\mu(\sigma\eta)=0.$
The first of these requires $\sigma^\vare\simeq\sigma.$  The second of these requires 
$(\sigma\eta)^\vare\simeq\sigma\eta.$
But since $\eta^\vare=\eta,$ we will again have $\sigma^\vare\simeq\sigma$ in this case.

Now suppose that $m>0$ and $\pi=\pi_1\otimes\tau.$ Then $\pi^w\simeq\pi$ implies $\pi_1^\vare\simeq\pi_1.$ Consider
$\pi$ and $\pi\chi_\eta.$ Then $\mu(\pi)=0$ if and only if $L(s,\pi_1, A)L(2s,\pi_1\times\tau)$ has a pole at $s=0$
(see \cite{\bf Sh2,GS}).
Similarly $\mu(\pi\chi_\eta)=0$ if and only if $L(s,\pi_1\chi_\eta, A)L(2s,\pi_1\eta\times\tau\eta)$ has a pole.
As in the case $m=0,$ precisely one of $L(s,\pi_1,A)$ or $L(s,\pi_1\chi_\eta,A)$ has a pole at $s=0,$ and thus, at least
one of $\mu(\pi)$ and $\mu(\pi\chi_\eta)$ is zero. But, as $\sigma$ is a component of both $\pi$ and $\pi\chi_\eta,$ we
must have $\mu(\sigma)=0,$  which therefore requires $\sigma^w\simeq\sigma.$
\qed\enddemo

\proclaim{Lemma 2.5}  Suppose $m=1.$  Then, for any $\pi\in\Cal E_c(\wt M)$ the representation $\sigma=\pi|_{M}$ is
irreducible.  Hence when $\pi^w\simeq\pi\chi,$ we have $\sigma^w\simeq\sigma.$
\endproclaim

\demo{Proof}  Since $m=1,$ we have
$$\wt M\simeq GL_{n_1}(E)\times GL_{n_2}(E)\times\dots\times GL_{n_r}(E)\times E^1.$$
Note that if $(g_1,g_2,\dots,g_r,h)\in M,$ then 
$$h=\left(\prod_{i=1}^r\det(g_i\vare(g_i))^{-1}\right).$$
Thus, $$M\simeq GL_{n_1}(E)\times\dots\times GL_{n_r}(E).$$
Now it is clear that if $\pi\simeq\pi_1\otimes\dots\otimes\pi_r\otimes\xi\in \var_c(\wt M),$
then $\ds{\pi|_{M}\simeq\left(\pi_1\otimes\dots\otimes\pi_r\right)\xi\xi^\vare}$ is irreducible.\qed\enddemo




\proclaim{Lemma 2.6}Suppose $m=0,$ or $\pi\simeq\pi_1\otimes\pi_2\otimes\dots\otimes\pi_r\otimes\tau,$ with
$\tau|_{SU(m)}$ of multiplicity 1.
\roster
\item"(a)"If $s\in \Cal S$ satisfies $\pi^s\simeq\pi\chi$
for some
$\chi\in X(\tM)$, then
$\sigma^s\simeq\sigma$ for any irreducible component $\sigma$ of $\pi|_M$.
\item"(b)" Suppose $w=sc,$ where $s=s_1s_2\dots s_k$ is the disjoint cycle decomposition, and $c$ changes an even
number of signs in each cycle $s_i.$  If $\pi^w\simeq\pi\chi,$ then $\sigma^w\simeq\sigma,$ for each component
$\sigma$ of $\pi|_M.$
\endroster
\endproclaim

\demo{Proof} 
(a) The argument is essentially that of Lemma 2.3 of \cite{\bf G1}.   
We give the proof in the case $m=0,$ and the
proof when $\rho|_{SU(m)}$ is multiplicity one is identical.
Let 
$$M_0=[M,M]=[\tM,\tM]\simeq SL_{n_1}(E)\times SL_{n_2}(E)\times\dots\times SL_{n_r}(E).$$
Since each $\pi_i|_{SL_{n_i}(E)}$ is multiplicity one \cite{\bf T}, so is $\pi|_{M_0}.$
If $\rho$ is a component of $\sigma|_{M_0},$ then $\rho\simeq\rho_1\otimes\rho_2\otimes\dots\otimes\rho_r,$ for some
choice of components $\rho_i$ of
$\pi_i|_{SL_{n_i}(E)}.$ 
Suppose $s=s_1s_2\dots s_k$ is the disjoint cycle decomposition of $s$ and, without loss of generality, assume that
$s_1=(1\,2\,\dots\,j_1)$ $s_2=(j_1+1\,\dots\,j_2),\dots,$ and $s_k=(j_{k-1}+1,\dots j_k).$  Let $j=j_1.$ Since
$\pi^s\simeq\pi\chi,$ we then have $\pi_{i+1}\simeq\pi_i\chi\simeq\pi_1\chi^i,$ for $i=1,2,...j-1,$ and $\pi_1\simeq\pi_j\chi,$
i.e., $\pi_1\simeq\pi_1\chi^j.$ Thus, for each $1\leq i\leq j,$  $\rho_i$ is an irreducible component of
$\pi_1.$ By \cite{\bf G-K}, for each $1\leq i\leq j-1,$ there is an $a_i\in E^\times,$ so that
$\rho_{i+1}=\rho_{i}^{\delta(a_i)},$ where
$\delta(a)=\pmatrix a\\&I_{n_1-1}\endpmatrix.$ Let $a_j=(a_1a_2\dots a_{j-1})^{-1}.$  Then $\rho_j^{\delta(a_j)}=\rho_1.$
Set $g_1=\diag\{\delta(a_1),\delta(a_2),\dots,\delta(a_j)\}.$  Then $\det g_1=1,$ and
$(\rho_1\otimes\dots\otimes\rho_j)^{g_1}=\rho_2\otimes\dots\otimes\rho_j\otimes\rho_1=(\rho_1\otimes\dots\otimes\rho_j)^{s_1}.$
Similarly, for $i=2,3,...k,$ we can find such a $g_i,$ with determinant 1. Setting
$g=\diag\{g_1,g_2,\dots g_k,\vare(g_k),\dots,\vare(g_1)\},$ we have $g\in M$ and $\rho\simeq\rho^g\simeq\rho^s.$  Therefore,  $\rho^s$ is a
component of both $\sigma|_{M_0}$ and $\sigma^s|_{M_0},$ and thus by multiplicity one, $\sigma^s\simeq\sigma.$

(b) Now suppose $w=sc,$ with $c\neq1.$  We again assume $s=s_1s_2\dots s_k.$  Suppose $s_1=(1\,2\,\dots\,j),$ and for
some
$d$ with
$0\leq d\leq j-1,$  that $c=C_{d+1}C_{d+2}\dots C_jc',$ where $c'$ acts trivially on $\{1,2,\dots,j\}.$
Then,
$\pi^w\simeq\pi\chi$ implies
$\pi_{i+1}^{\vare^{b_i}}\simeq\pi_i\chi$
for $i=1,2,\dots, j-1,$ and $\pi_1^{\vare^{b_j}}\simeq\pi_j\chi.$ Here $b_i\in\{0,1\}.$
Let $\rho=\rho_1\otimes\dots\rho_r\otimes\tau,$ be a component of $\pi|_{M_0},$ where, again, $M_0=[\wt M,\wt M].$  Then, for
$i=1,\dots,j,$ we have
$(\rho_{i+1})^{\vare^{b_i}}=\rho_i^{\delta(a_{i})},$
for some $a_i.$  
Since $j-d$ is even, we let 
$$b=\vare\left(\delta\left(a_1\right)^{-1}\delta(a_2)^{-1}\dots\delta(a_d)^{-1}\right)(\delta(a_{d+1})^{-1}\vare(\delta(a_{d+2}))^{-1}
(\delta(a_{d+3}))^{-1}\dots\vare(\delta(a_{j-1}))^{-1}.$$ Then $\rho_j^b\simeq\rho_1^{\vare}.$ 
Therefore, we set
$$g_1=\diag\{\delta(a_1),\dots,\delta(a_{j-1}),b\}.$$
Note that $\det(g_1)\det(\vare(g_1))=1,$ and thus, choosing $g_2,\dots,g_k$ in a similar manner, we have
$$g=\diag\{g_1,g_2,\dots,g_k,\vare(g_k),\dots,\vare(g_1)\}\in M,$$
with $\rho^g\simeq\rho^w.$ Therefore, we again see that $\rho^w\simeq\rho.$\qed
\enddemo

\proclaim{Lemma 2.7}  Suppose $\pi=\pi_1\otimes\pi_2\otimes\dots\otimes\pi_r\otimes\tau,$ and
that $\tau|_{SU(m)}$ is multiplicity one.  If $c\in\Cal Z$ satisfies $\pi^c\simeq\pi,$ then $\sigma^c\simeq\sigma,$ for
each component $\sigma$ of $\pi|_{M}.$
\endproclaim

\demo{Proof}  Let $M_0=M_1\times M_2\times\dots\times M_r\times SU(m)\subset M,$
with $M_i=\{g\in GL_{n_i}(E)|\det g\in F^\times\}.$  Then, $\pi|_{M_0}$ is multiplicity one.
If $\ds{c=\prod_{i\in B}C_i,}$ then
$\ds{\pi^c\simeq\bigotimes_{i\in B}\pi_i^\vare\otimes\bigotimes_{i\not\in B}\pi_i\otimes\tau.}$
Let $\sigma$ be a component of $\pi|_{M},$ and suppose $\sigma_0$ is a component of $\sigma|_{M_0}.$
Then $\sigma_0=\rho_1\otimes\dots\otimes\rho_r\otimes\tau_0,$ where $\rho_i$ is a component of $\pi_i|_{M_i}$ and $\tau_0$ is
a component of $\tau|_{SU(m)}.$  Therefore,
$$\sigma_0^c\simeq\bigotimes_{i\in B}\rho_i^\vare\otimes\bigotimes_{i\not\in B} \rho_i\otimes\tau_0\simeq\sigma_0,$$
by Lemma 2.4.  Thus, $\sigma_0$ is a component of both $\sigma$ and $\sigma^c$ upon restricition to $M_0,$ and hence by
multiplicity one $\sigma\simeq\sigma^c.$\qed\enddemo

\remark{Remark} If $\tau$ is generic, then the representation $\pi$ satisfies the hypotheses of Lemmas 2.6 and 2.7.
\endremark

Now assume $m\geq 2$. Recall
that, for
$h\in U_m(F)$, we have det
$h\in E^1$, and thus, by Hilbert's Theorem 90, det $h=a\overline a^{-1}$, for some $a\in E^\times.$ For $a\in E^\times$ we
let
$$
\alpha_m(a)=\bmatrix a\\ &I_{m-2}\\ &&\overline a^{-1}\endbmatrix.
$$ 
Then $\det(\alpha_m(a))=\det h$. Note $\alpha_m(ab)=\alpha_m(ba)=\alpha_m(a)\alpha_m(b).$
If $g\in GL_k(E)$, we abuse notation and write $\alpha_m(g)$ for $\alpha_m(\det(g))$.

\proclaim{Lemma 2.8}If $\widetilde M\cong GL_{n_1}(E)\times GL_{n_2}(E)\times\ldots\times GL_{n_r}(E)\times U_m(F)$, and
$m\geq 2$, then,
$M\simeq GL_{n_1}(E)\times\ldots\times GL_{n_r} (E)\rtimes SU_m(F)$.
\endproclaim

\demo{Proof}Let $g=(g_1,g_2,\ldots,g_r)\in GL_{n_1}(E)\times\ldots\times GL_{n_r}(E)$. Then the map
$(g,h_0)\mapsto\pmatrix g\\ &\alpha_m(g)^{-1}h_0\\ &&\vare(g)\endpmatrix$, is a set bijection from
$GL_{n_1}\times\ldots\times GL_n(E)\rtimes SU_m(F)$ to $M$. The action of $GL_{n_1}(E)\times\ldots\times GL_{n_r}(E)$ on
$SU_m(F)$ is by $g\circ h=\alpha_m(g) h\alpha_m(g)^{-1}$. It is now clear to see our map is an isomorphism.\qed
\enddemo

If $x=(g_1,\dots,g_r,h)=(g,h)\in\tilde M,$ and $w\in W,$ then
we denote $wxw^{-1}$ by $(g^w,h).$  Restricting to $M,$ we see that, under the above isomorphism,
$W$ acts on $(GL_{n_1}(E)\times\dots\times GL_{n_r}(E))\ltimes SU_m(F)$ by
$w\cdot(g,h_0)=(g^w,(\alpha_m(g)^{-1}\alpha_m(g^w)h_0).$  (Note that $\alpha_m(g^{-1})\alpha(g^w)\in
SU_m(F).$)

Now suppose $\pi\simeq\pi_1\otimes\pi_2\otimes\pi_2\dots\otimes\pi_r\otimes\tau$ and that
$V=V_1\otimes V_2\otimes\dots\otimes V_r\otimes V_{\tau}$ is the space of $\pi.$  We may write $g=(g_1,\dots,g_r)\in
GL_{n_1}(E)\times\dots\times GL_{n_r}(E),$  and $\pi_0(g)=\otimes_i\pi_i(g_i)$ acting on $\ds{V_0'=\otimes_{i=1}^rV_i.}$  Let
$(\tau_0, V_0)$ be a component of
$\tau|_{SU_m(F)}.$ Then,  with respect to the semidirect product dercomposition in Lemma 2.8, 
the map
$$(g,h_0)\mapsto\left( v_0'\otimes v_0\mapsto\pi_0(g)v_0'\otimes\tau(\a_m(g)^{-1})\tau_0(h_0) v_0\right),\tag
2.6$$ is an irreducible component of $\pi|_M.$  We now prove another partial converse to Corollary 2.2.  Note that 
we do not  assume multiplicity one upon restriction.  This result is crucial to the $R$--group computations of Section 3.

\proclaim{Proposition 2.9}  
Suppose $M\simeq GL_{n_1}(E)\times\dots\times GL_{n_r}(E)\times U_m(F)$ and $m\geq 2.$
Let $\pi\in \Cal E_2(\tilde M)$
and suppose
$\sigma$ is an irreducible component of
$\pi|_M.$ Then $W(\pi)\subset W(\sigma).$
\endproclaim

\demo{Proof} Let $w\in W(\pi)$ and suppose that $w=sc,$ with $s\in\Cal S$ and $c\in\Cal Z.$
We note that if $\pi\simeq\pi_1\otimes\dots\otimes\pi_r\otimes\tau,$ then
$$\pi^w\simeq\pi_{s(1)}^{\vare_1}\otimes\pi_{s(2)}^{\vare_2}\otimes\dots\otimes\pi_{s(r)}^{\vare_r}\otimes\tau,$$
where each $\vare_i$ is either $\vare$ or trivial.
Since $\pi^w\simeq\pi,$ we have $\pi_i\simeq\pi_{s(i)}^{\vare_i}$ for each $i,$ so fix an intertwining,
$T:V_{s(i)}\rightarrow V_i$ with $\pi_i T_i=T_i\pi_{s(i)}^{\vare_i}.$  Then $\ds{T=\bigotimes_i T_i\otimes
\text{id}_{V_\tau}}$ satisfies $\pi T =T\pi^w.$  Now suppose that 
$\sigma$ is given by (2.6).
Then 
$$\gather
T\sigma^w(g,h_0)(v_0'\otimes v_0)=\bigotimes_{i=1}^r T_i\pi_{s(i)}^{\vare_i}(g_i)v_i\,\otimes
\tau(\a_m(g^w)^{-1})\tau_0(\alpha_m(g)^{-1}\alpha_m(g^w)h_0) v_0\\=\bigotimes_{i=1}^r \pi_i(g_i) T_i(v_i)\,\otimes
\tau(\a_m(g)^{-1})\tau_0(h_0) v_0 =\sigma(g,h_0)T(v_0'\otimes v_0).
\endgather$$
Thus $T$ is an equivalence between $\sigma^w$ and $\sigma.$\qed\enddemo

\remark{Remarks} \item{(a)} The proof of Proposition 2.9 will apply to to the case $\pi^w\simeq\pi\chi$ whenever
the identity is an equivalence between $\tau$ and $\tau\chi.$

\item{(b)} There is a similar result when $m=0,$ using an identification $M\simeq (GL_{n_1}(E)\times
GL_{n_2}(E)\times\dots\times GL_{n_{r-1}}(E))\ltimes SL_{n_r}(E).$  To the extent we need this result in Section 3, however,
Lemma 2.7 will suffice.
\endremark

\proclaim{Corollary 2.10} For any $\pi\in\Cal E_2(\tilde M),$ and any irreducible component $\sigma$ of $\pi|_M,$ we
have $W(\pi)\subset W(\sigma).$
\endproclaim

\demo{Proof}  This will follow immediately from Lemma 2.4, 2.6, 2.7, and Proposition 2.9, unless
$m=0$ and $w$ does not satisfy the hypotheses of either part of Lemma 2.6.  This case we now resolve.
All that is left is to consider the case where $m=0$ and $w=sc,$ with $c\neq1$ changing an odd number of signs in some cycle of $s.$
So, suppose $s=s_1s_2\dots s_k$ is the disjoint cycle decomposition for $s,$ with $s_1=(12\dots j).$  Further suppose that
for some $1\leq d\leq j,$ we have $c=C_dC_{d+1}\dots C_jc',$ with $c'$ acting trivially on $\{1,\dots,j\}$ and $j-d+1$ odd.
Then $\pi^w\simeq\pi$ implies
$$\pi_1\simeq\pi_2\simeq\dots\simeq\pi_d\simeq\pi_{d+1}^\vare\simeq\pi_{d+2}\simeq\dots\simeq\pi_j\simeq\pi_1^\vare.$$
Therefore, $\pi_i\simeq\pi_i^\vare$ for $1\leq i\leq j.$  Hence, $C_j\in W(\pi).$  Taking $w_1=wC_j,$ we see that $w_1$
changes an even number of signs among $\{1,\dots,j\}.$  Proceeding in the same manner on the cycles $s_2,\dots,s_k,$ if necessary,
we see we can write $w=w_0c_0,$ with $w_0$ satisfying the hypotheses of Lemma 2.6(b) and both $w_0$ and $c_0$ in $W(\pi).$
Since $m=0,$ Lemma
2.7 applies to
$c_0.$  Thus, we know
$\sigma^{w_0}\simeq\sigma$ and $\sigma^{c_0}\simeq\sigma,$ and hence $\sigma^w\simeq\sigma.$\qed\enddemo

\noindent {\bf Section 3.\ $R$--groups}

For the moment we work in the more general setting where $G\subset \wt G$ have the same derived group. 
Let $R_{\pi}(\sigma)=R(\pi)\cap W(\sigma).$ Note that if
$w\in R_{\pi}(\sigma), $ then, $w\Delta'(\pi)=\Delta'(\pi),$ 
and thus by Lemma 2.3, $w\Delta'(\sigma)=\Delta'(\sigma),$
which, combined with the fact that $w\in W(\sigma),$ shows that $w\in R(\sigma).$ Note that if $w\in R(\sigma)\cap
W(\pi),$ then, as $w\Delta'(\sigma)=\Delta'(\sigma),$ we have $w\in R_\pi(\sigma),$ so we certainly have
$R_\pi(\sigma)=R(\sigma)\cap W(\pi).$
In fact
$R_\pi(\sigma)$ is a normal subgroup of
$R(\sigma).$  In this section we first describe the quotient $R(\sigma)/R_\pi(\sigma).$  Then, specializing to the
case where $\tilde G=U_n,$
and  $G=SU_n,$  we show that $R(\sigma)$ is in fact a semidirect product of $R(\pi)$ and a naturally occuring group of
characters.

\definition{Definition 3.1 }Let $\widehat{W(\sigma)}=\{\chi|w\pi\simeq\pi\chi$ for some $w\in W(\sigma)\}$. 
\enddefinition

Let $\chi\in\widehat{W(\sigma)}/X(\pi)$. Choose $w\in W(\sigma)$ with $\pi^w\simeq\pi\chi.$  Then $w=rw',$ with $r\in
R(\sigma)$ and $w'\in W'(\sigma).$ By Lemma 2.3 $w'\in W'(\pi), $ and in particular,  $\pi^{w'}\simeq\pi,$  we must have
$\pi^r\simeq\pi\chi.$   Hence, for any $\chi\in\widehat{W(\sigma)}$ there is an element $r\in R(\sigma)$ with
$\pi^r\simeq\pi\chi.$
The following result is now obvious.

\proclaim{Proposition 3.2}  For any $\pi\in\var_2(\wt M)$ and any irreducible component $\sigma$ of $\pi|_{M}$ we have
$R_{\pi}(\sigma)\triangleleft R(\sigma)$ and $R(\sigma)/R_\pi(\sigma)\simeq \widehat{W(\sigma)}/X(\pi).$\qed
\endproclaim

We now specialize tot he case where $\wt G=U_n(F)$ and $G=SU_n(F).$  For $C\in\cZ$
there is a
$B\subseteq \{1,2,\ldots r\}$ with $C=C_B=\prod\limits_{i\in B} C_i$.   By \cite{\bf G2} there is a subset
$B(\pi)$ so that $R(\pi)=\left<C_i|i\in B(\pi)\right>.$    By Corollary 2.10, $R(\pi)\subset R(\sigma).$  Also note that
if $\pi^{w_1}\simeq\pi^{w_2}\simeq\pi\chi,$ then $\pi^{w_1^{-1}w_2}\simeq\pi,$ 
which implies $\sigma^{w_1^{-1}w_2}\simeq\sigma.$
Hence if $\chi\in \widehat{W(\sigma)},$
then $\sigma^w\simeq\sigma$ for all $w\in W$ with $\pi^w\simeq\pi\chi.$

\proclaim{Lemma 3.3}
Suppose $\chi\in\widehat{W(\sigma)}$ and further suppose there are
elements
$w_1=s_1 c_1,\ w_2=s_2 c_2\in R(\sigma)$ with $\pi^{w_1}\simeq\pi\chi\simeq
\pi^{w_2}$. Here
$s_i\in\cS$, and
$c_i\in\cZ$. Then $s_1=s_2$.
\endproclaim

\demo{Proof}
Note $w_1 w_2^{-1}\in R(\sigma)\cap W(\pi)=R(\pi).$   Therefore, $w_2w_2^{-1}\in\Cal
Z.$ Since
$w_1 w_2^{-1}=(s_1 s_2^{-1})(s_2 c_1 s_2^{-1} c_1)$ is the decomposition of $w_1 w_2^{-1}$ in the semidirect product
$W=\cS\ltimes\cZ$, we have $s_1=s_2$.\qed\enddemo

For $\chi\in \widehat{W(\sigma)},$ denote by $s_\chi$ the unique element in $\cS$
so that $s_\chi c\in R(\sigma)$ for some $c\in\cZ.$   
We will give an explicit description of $s_\chi.$  First, we need a lemma.

\proclaim{Lemma 3.4} Suppose $\pi\simeq\pi_1\otimes\dots\otimes\pi_r\otimes\tau.$  Let $w=sc\in W(\bold G,\bold
A)$ with $\pi^w\simeq\pi\chi,$ for some $\chi.$  Suppose, $c(i)\neq i.$  Then $\pi_i\chi\simeq\pi_i\chi^\vare.$
\endproclaim

\demo{Proof} 
First suppose that $s(i)=i,$ so $w(i)=-i.$ Then $\pi^w\simeq\pi\chi$ implies $\pi_i^\vare\simeq\pi_i\chi.$  Applying
$w$ again we see $\pi_i\simeq\pi_i\chi^2,$ while $\pi_i\simeq
(\pi_i\chi)^\vare\simeq\pi_i^\vare\chi^\vare\simeq\pi_i\chi\chi^\vare.$ Thus, $\pi_i\chi\simeq\pi_i\chi^\vare.$

Now suppose that $s(i)\neq i.$  Assume
$s=s_1s_2\dots s_\ell,$ is the disjoint cycle
decomposition of $s,$ and $s_1(i)\neq i.$ Without loss of generality, suppose that $s_1=(12\dots j),$ and that, for some
$1\leq d\leq j,$ we have
$c=C_dC_{d+1}\dots C_jc',$ with
$c'$ trivial on $\{1,\dots,j\}.$
If $d\geq 2,$ then $\pi^w\simeq\pi\chi$ implies $\pi_2\simeq\pi_1\chi.$  Then $\pi^{w^2}\simeq\pi\chi^2,$ implies
$\pi_2^\vare\simeq\pi_j^\vare\chi^2\simeq\pi_1^\vare\chi.$  Thus $\pi_2\simeq\pi_1\chi^\vare\simeq\pi_1\chi,$ proving
the claim in this case.  If $d=1,$ then $\pi^w\simeq\pi\chi$ implies $\pi_2^\vare\simeq\pi_1\chi,$ while applying $w^2$
gives
$\pi_2\simeq\pi_1^\vare\chi\simeq\pi_1^\vare\chi^\vare,$ proving the claim in this case as well.\qed\enddemo

We now describe the permutation $s_\chi$ explicitly. For $\chi\in \widehat{W(\sigma)},$ we let
$\Omega(\chi,\pi)=\{i|\pi_i\simeq\pi_i\chi,\text{ or }\pi_i^\vare\simeq\pi_i\chi\}.$ We then let
$\O_1(\chi,\pi)=\{1,2,\dots,r\}\setminus\O(\pi,\chi).$
For $i\in\O(\chi,\pi),$ let $s(i)=i.$ If
$\O_1(\chi,\pi)=\emptyset,$ we are done and $s=1.$  Otherwise, for each $i\in\O_1(\chi,\pi),$ we let
$\O_1(i,\chi,\pi)=\{j\in\O_1(\chi,\pi)|\pi_j\simeq\pi_i\chi\},$ and let
$\O_1^\vare(i,\chi,\pi)=\{j\in\O_1(\chi,\pi)|\pi_j^\vare\simeq\pi_i\chi\}.$ Let $i_{11}=\min(\O_1(\chi,\pi)).$
Define 
$$i_{12}=\cases\min(\O_1(i_{11},\chi,\pi))&\text{if
}\O_1(i_{11},\chi,\pi)\neq\emptyset,\\\max(\O_1^\vare(i_{11},\chi,\pi))&\text{otherwise.}\endcases$$
Suppose we have defined $i_{11},i_{12},\dots i_{1j},$ if $i_{1j}=i_{11},$ then we let
$s_1=(i_{11},i_{12},\dots,i_{1j-1}).$ Otherwise, we let
$$i_{1j+1}=\cases\min(\O_1(i_{1j},\chi,\pi))&\text{if
}\O_1(i_{1j},\chi,\pi)\neq\emptyset,\\\max(\O_1^\vare(i_{1j},\chi,\pi))&\text{otherwise.}\endcases$$
This, inductively, defines an element $s_1$ of $\Cal S$

Now let $\O_2(\chi,\pi)=\O_1(\chi,\pi)\setminus\{i_{11},\dots,i_{1j-1}\}.$ If $\O_2(\chi,\pi)=\emptyset,$ we are done,
and $s=s_1.$  Otherwise, repeat the above process by taking $i_{21}=\min(\O_2(\chi,\pi)),$ and defining
$$i_{2j+1}=\cases\min(\O_1(i_{2j},\chi,\pi))&\text{if
}\O_1(i_{2j},\chi,\pi)\neq\emptyset,\\\max(\O_1^\vare(i_{2j},\chi,\pi))&\text{otherwise.}\endcases$$
and let $j_2$ be the smallest integer greater than $1$ for which $i_{2j_2}=i_{21}.$  Set
$s_2=(i_{21},\dots,i_{2j_2-1}).$ 
Proceed inductively to define $\O_3(\chi,\pi), \O_4(\chi,\pi),\dots \O_k(\pi,\chi),$ with associated cycles
$s_3,s_4,\dots s_k,$ and suppose that $k$ is minimal with the property that $\O_{k+1}(\pi,\chi)=\emptyset.$
Then let $s=s_1s_2\dots s_k.$  By construction, if $\pi_i\simeq\pi_j$ for $i<j,$ and $s_{i'}(i)\neq i,$
then $s_{j'}(j)\neq j$ for some $j'>i'.$

\proclaim{Lemma 3.5}  Let $\chi\in \widehat{W(\sigma)},$ and define $s=s_1s_2\dots s_k,$ as above. Then $s=s_\chi.$
\endproclaim

\demo{Proof}
We define $w=sc$ with the property that $\pi^w\simeq\pi\chi,$ and then show $w\in R(\sigma).$ That will show that
$s=s_\chi.$  We let $c$ be defined by 
$$c(i)=\cases i&\text{if }\pi_{s(i)}\simeq\pi_i\chi\\-i&\text{if }\pi_{s(i)}^{\vare}\simeq
\pi_i\chi\text{ and }\pi_{s(i)}\not\simeq\pi_i\chi.\endcases$$ Then, setting $w=sc,$ we have $\pi^w\simeq\pi\chi$
by construction.  We show that
$w$ preserves the positivity of all elements of $\Delta'(\sigma).$  First suppose
$\alpha=\alpha_{ij}=e_{b_i}-e_{b_{j-1}+1},$ where $b_i=n_1+n_2+\dots n_i.$  Then, $\pi_i\simeq\pi_j.$  If $w(i)=s(i),$
then $w\alpha_{ij}=e_{s(i)}\pm e_{s(j)}.$  This would be a negative root if and only if $w(j)=s(j)$ and $s(j)<s(i).$
But as $w(i)=s(i),$ we know if $w(j)=s(j),$ then $s(j)>s(i),$ by construction.  On the other hand, if $w(i)=-s(i),$
then we must have
$s(j)<s(i)$ and $c(j)=-j,$ as well, so again $w\alpha>0.$

Suppose that $\alpha=\alpha'_{ij}=e_{b_i}+e_{b_{j-1}+1}\in\Delta'(\sigma).$  Then
$\pi_j\simeq\pi_i^\vare.$ If $w\alpha_{ij}'<0,$ then certainly either $w(i)=-s(i),$ or $w(j)=-s(j).$  Thus, by
Lemma 3.4,
$\pi_i\chi\simeq\pi_i\chi^\vare.$  If $w(i)=s(i),$ and $w(j)=-s(j),$ then
$\pi_{s(i)}\simeq\pi_i\chi.$  Note that $\pi_{s(j)}\simeq\pi_j^\vare\chi^\vare\simeq\pi_i\chi,$ and thus, by the
choice of $s(i), $ we must have $s(i)<s(j).$  Therefore $w\alpha'_{ij}>0,$ contradicting our assumption.  So we must
have $w(i)=-s(i).$  If $w(j)=-s(j),$ then $\pi_{s(j)}^\vare\simeq\pi_j\chi\simeq\pi_i^\vare\chi$ which implies
$\pi_{s(j)}\simeq\pi_i\chi.$ Then, by construction of $s,$ we would have assigned $s(i)=j',$ with $j'\leq s(j),$ and
$\pi_{j'}\simeq\pi_i\chi.$  We would then have $w(i)=s(i)=j',$ contradicting our assumption. 
Thus, we must have
$w(j)=s(j).$  Then
$\pi_{s(j)}\simeq\pi_j\chi\simeq\pi_i^\vare\chi\simeq(\pi_i\chi)^\vare,$ and thus by construction, $s(i)>s(j).$  Therefore,
$w\alpha_{ij}'=\alpha_{s(j)s(i)}>0,$ which contradicts our assumption.  Hence it is impossible for
$w{\alpha_{ij}'}<0$ if $\alpha_{ij}'\in\Delta'(\sigma).$ 

Finally, suppose that $\ds{\beta_i=\cases e_{b_i} &\text{if } G=U(2k+1)\\2e_{b_i}&\text{if }G=U(2k).\endcases}$
Then $\pi_i^\vare\simeq\pi_i.$  Assume $w\beta_i<0.$   Then we must have $c(i)=-i.$ If $s(i)=i,$ then
$\pi_i^\vare\simeq\pi\chi,$ which now says
$\pi_i\simeq\pi_i\chi,$ and by assumption we have $c(i)=i,$ contradicting our assumption.  Now suppose that $s(i)\neq
i.$ Then as $c(i)=-i,$ we have, by Lemma 3.4, $\pi_i\chi\simeq\pi_i\chi^\vare,$ and  
$\pi_{s(i)}^\vare\simeq\pi_i\chi.$
Since
$\pi_i^\vare\simeq\pi_i,$ we have
$\pi_{s(i)}\simeq\pi_{i}\chi,$ which means that $\O_{\ell}(i,\chi,\pi)\neq\emptyset,$ and therefore,  we are forced to
take $c(i)=i.$  This also contradicts our assumption, and hence $w\beta_i>0$ for all
$\beta_i\in \Delta'(\sigma).$

In order to conclude that $s=s_\chi,$ we must know that $w\sigma\simeq\sigma.$ 
However, by assumption, $\chi\in\widehat{W(\sigma)},$ and thus $\sigma^w\simeq\sigma.$ (See the remark preceding Lemma
3.3.)\qed\enddemo

\remark{Remark}  Note that it is possible that $s_\chi=s_\eta$ with $\pi\eta\not\simeq\pi\chi.$  In particular, it is
possible that $s_\chi=1,$ and $c\in\Cal Z$ is an element of $W(\sigma)\setminus W(\pi).$
\endremark

\proclaim{Lemma 3.6}  Suppose  $\chi\in \widehat{W(\sigma)}/X(\pi).$
Then there is a unique
minimal
$B_\chi$ with
$s_\chi C_{B_\chi}\in R(\sigma).$
\endproclaim

\demo{Proof}
Choose any element $w\in R(\sigma)$ with $w\pi\simeq\pi\chi$.  By Lemma 3.3, we have $w=s_{\chi}C_B$ for some $B.$ 
Let
$B'=\{i\in B|C_i\in R(\pi)\}$. Let $B_\chi=B\backslash B'$. Then,   $C_{B'}\in R(\sigma)$,
and thus $w_\chi=s_\chi C_{B_\chi}=s_\chi C_B C_{B'}\in R(\sigma)$. For any $B_1$, with $s_\chi C_{B_1}\in R(\sigma),$
we have
$(s_\chi C_{B_\chi})^{-1} s_{\chi}C_{B_1}=C_{B_2}\in R(\pi)$ so $B_2\subset B(\pi),$ and 
$s_\chi C_{B_1}=s_{\chi}C_{B_\chi}C_{B_2}.$ Therefore, as $B_\chi\cap B(\pi)=\emptyset,$
$B_\chi\subset B_1.$\qed 
\enddemo

For $\chi\in \widehat{W(\sigma)}$ we let $w_\chi$ be the element of $R(\sigma)$ given by Lemma 3.4.
We now show that the sequence $1\rightarrow R(\pi)\rightarrow R(\sigma)\rightarrow \widehat{W(\sigma)}\rightarrow 1$
splits, so that $R(\sigma)$ is a semidirect product.

\proclaim{Theorem 3.7}Let $\chi_1,\chi_2\in\widehat{W(\sigma)}/X(\pi).$  
Then $w_{\chi_1} w_{\chi_2}=w_{\chi_1\chi_2}$.  Thus $\Gamma_\sigma=\left\{w_\chi\,|\,\chi\in
\widehat{W(\sigma)}\right\}$ is a subgroup of $R(\sigma)$ and $R(\sigma)=\Gamma_\sigma\ltimes R(\pi).$
\endproclaim

\demo{Proof}  Let $w_i=w_{\chi_i}$
and suppose $w_i=s_i C_{B_i}$ is the
decomposition given in Lemma 3.6. Then $B_i\cap B(\pi)=\emptyset$. Note that $w_1 w_2=s_1 s_2
C_{s_2^{-1}(B_1)}C_{B_2}$. Suppose
$s_2^{-1} (B_1)\cap B(\pi)\not=\emptyset$. Let $j\in s_2^{-1}(B_1)\cap B(\pi).$
Then  $C_{s_2^{-1}(B_1)}C_j$ has shorter length than
$C_{s_2^{-1}(B_1)}$, hence shorter length than $|B_1|$. Let $C_{s_2^{-1}(B_1)}C_j=C_{B'_1}$. Then
$$ w'=w_1 w_2 C_j w_2^{-1}=s_1 s_2 C_j C_{s_2^{-1}(B_1)}C_{B_2} C_{B_2} s_2^{-1}=s_1 C_{s_2(B_1')},
$$ and $\pi^{w'}\simeq\pi\chi$. However, since $|s_2(B'_1)|< |B_1|$, this contradicts our choice of $B_1$. Thus,
$s_2^{-1}(B_1)\cap B(\pi)=\emptyset$. Since $B_2\cap B(\pi)=\emptyset$, as well we see $C_{s_2^{-1}(B_1)}
C_{B_2}=C_B$, with $B\cap B(\pi)=\emptyset$. Now
$w_1 w_2=s_1 s_2 C_B=s_{\chi_1\chi_2}C_B\in R(\sigma)$, and since $B\cap B(\pi)=\emptyset$, we know $B$ is minimal with respect to this
property. Thus, $w_1 w_2=w_{\chi_1\chi_2}$.\qed
\enddemo

We address non-abelian $R$--groups. 
Suppose
$w_\chi=s_\chi C_{B_\chi}\in
\Gamma_{\sigma}.$ Suppose  further that $s_{\chi}(i)=i$ for all $i\in B(\pi).$  Then
$w_\chi$ centralizes $R(\pi),$ and hence lies in the center of $R(\sigma).$  Thus, $R$ is non-abelian if and only
if there is some $w=s_\chi c\in R(\sigma)$ with $s(i)\neq i$ for some $i\in B(\pi).$

\proclaim{Lemma 3.8}
Let $R(\pi)=\,\left<C_i|i\in B(\pi)\right>.$  
\roster
\item"(a)" There is no $j\not\in B(\pi)$ for which $\pi_j^\vare\simeq\pi_j$ and $C_j\in R(\sigma).$
\item"(b)" If $w=sc\in R(\sigma)$ then $s(B(\pi))=B(\pi).$
\endroster
\endproclaim
\demo{Proof}
(a) Suppose $j\notin B(\pi)$ and $\pi_j^\vare\simeq\pi_j.$ Then $C_{j}\in W(\pi)$ and $C_j\not\in R(\pi),$ by
assumption.  Thus, $C_j\alpha<0$ for some $\alpha\in \Delta'(\pi)=\Delta'(\sigma),$ which shows that $C_j\not\in
R(\sigma).$ 

For (b), suppose that there is some $i\in
B(\pi)$ with
 $s(i)\not\in B(\pi).$ Then, since $C_i\in R(\pi)\subset R(\sigma),$ we see that
$wC_i\in R(\sigma).$  Let $c=C_B.$  If $i\in B,$ then $wC_i=sC_{B\setminus\{i\}}\in R(\sigma),$ and thus we may assume
that
$i\not\in B.$ Now since $wC_i=C_{s(i)}w\in R(\sigma),$ we see that $C_{s(i)}\in R(\sigma),$ as well. 
Note that if $\pi^{C_{s(i)}}\simeq\pi\eta,$ then $\pi_{s(i)}^\vare\simeq\pi_{s(i)}\eta,$ while $\pi_i\simeq\pi_i\eta.$
Since $\pi_{s(i)}\simeq\pi_i\xi,$ for some $\xi,$ we have $\pi_{s(i)}\eta\simeq\pi_{s(i)}.$  Thus,
$\pi^{C_{s(i)}}\simeq\pi.$
However, this contradicts part (a), which completes the lemma.  
\qed
\enddemo

We now wish to describe the conditions under which $i_{G,M}(\sigma)$ has elliptic constituents when $\sigma$ is generic.
With this assumption, the cocycle is a coboundary (see \cite{{\bf K2}, page 62}), so $\Cal C(\sigma)\simeq\Bbb
C[R(\sigma)].$    We first, however, describe all the regular elements of $R(\sigma),$ without regard to whether
$\sigma$ is generic.  If we know that the cocycle split in general, this would then describe the elliptic spectrum
in general.

\proclaim{Theorem 3.9}  Let $\pi\in \Cal E_2(M),$ and suppose that $\sigma$ is a component of $\pi|_M.$
Then $R(\sigma)_{reg}\neq\emptyset$  if an only if there is an element $w=s_\chi C_B\in R(\sigma)$
with
$s_\chi$ an $r$-cycle and  $|B|$ odd.
\endproclaim

\demo{Proof}  If $s_{\chi}$ has more than one orbit, one can easily construct  a non-zero element $X\in \frak a_w.$
Let $O_1=\{1,s(1),s^2(1),\dots,s^{k-1}(1)\}$ and $O_2=\{j,s(j),\dots,s^{\ell-1}(j)\}$ be two distinct orbits of
$s_\chi.$
Let $a_i=0,$ for $i\not\in O_1\cup O_2.$ We set $a_1=a_j=1,$ and for $1\leq i\leq k=2,$ let
$$a_{s^i(1)}=\cases a_{s^{i-1}(1)}&\text{if } s^{i-1}(1)\not\in B;\\-a_{s^{i-1}(1)}&\text{if } s^{i-1}(1)\in
B.\endcases$$
similarly, for $1\leq i\leq \ell-2$ we let
$$a_{s^i(j)}=\cases a_{s^{i-1}(j)}&\text{ if } s^{i-1}(j)\not\in B;\\-a_{s^{i-1}(j)}&\text{if } s^{i-1}(j)\in
B.\endcases$$
Let $X=Y_M(a_1,a,\dots, a_r),$ as in Section 1.  Then $wX=X,$ and $tr X=0,$ so $X\in \frak a_w.$  Therefore, if $w=s_\chi C_B$ is
regular, then
$s_\chi$ is an $r$--cycle.

Suppose, without loss of generality, that $s_\chi=(12\dots r),$  
and $X=Y_M(a_1,\dots,a_r)\in \frak a_w.$  Suppose $a_1=\lambda.$  Then, for $2\leq i\leq r,$ we have $a_i=\pm\lambda.$
Thus, we simply denote $X=X_B(\lambda).$
Now, as $tr X=0,$ we have $nr(\lambda-\bar\lambda)=0,$ so $\lambda\in \Bbb R.$
 Note that if
$a_r=\lambda$ and $r\in B$ or $a_r=-\lambda$ and $r\not\in B,$  then $\lambda=-\lambda,$ and hence $\frak
a_w=\{0\}.$ Suppose $|B|$ is even. If $r\in B,$ then $a_r=-\lambda,$ as
$w$ changes an odd number of signs among
$1,2,\dots,r-1.$ Then $X_B(\lambda)\in\frak a_w$ for any $\lambda\in\Bbb R.$
Similarly, if $r\not\in B,$ then $a_r=\lambda$ and hence $X_B(\lambda)\in \frak a_w$ for any $\lambda\in\Bbb R.$  Thus,
$w\not\in R(\sigma)_{reg}$ if $|B|$ is even.
On the other hand, if $|B|$ is odd and $r\in B,$ then $a_r=\lambda,$ so $X_B(\lambda)\in \frak a_w$ only for $\lambda=0.$
Similarly, if $|B|$ is odd and $r\not\in B,$ then $ a_r=\lambda,$ so again $\frak a_w=\{0\}.$
\qed\enddemo

\remark{Remark 3.10} We note that Theorem 3.9 is inconsistent with the results of \cite{\bf G2}, and this is due to an error in that
work which we now correct.  Note that $\frak a_{\tilde G}=\{\lambda I_n|\lambda\in i\Bbb R\}.$  Then we easily see that 
$\frak a_w=\frak a_{\wt G}$ if and only if $m=0,$ $r=1,$ and $R\simeq\Bbb Z_2.$  That is, only the Siegel Levi subgroup of $U_n$
supports non-discrete elliptic representations.\qed\endremark

We note that by Lemma 3.6(b), there are two possibilities for $R(\sigma)$ in the situation of Theorem 3.9.  Namely, if
$R(\pi)=\{1\}$ then $s_\chi$ an $r$--cycle implies $R(\sigma)\simeq\Bbb Z^r.$  Otherwise $R(\pi)\neq\{1\}$
implies $R(\pi)=\,\left<C_i|\, 1\leq i\leq r\right>\simeq\Bbb Z_2^r.$ Now if $\eta\in \widehat{W(\sigma)},$ then $s_\eta\in
R(\sigma).$  If $s_{\eta}(1)\neq 1,$ and $s_{\eta}(j)=j,$ then $\pi_j\simeq\pi_j\eta$ implies $\pi_1\eta\simeq\pi_1,$
which contradicts our assumption that $s_{\eta}(1)\neq 1.$  This implies that $s_{\eta}=s_{\chi}^i$ for some $i,$ and
so  $R(\sigma)\simeq\Bbb Z_r\ltimes\Bbb Z_2^r.$

We now give a specific case of this second phenomenon, and note that the induced representation will contain both
elliptic and non-elliptic representations.  Suppose
$r=3$ and we are in this second case, namely
$R(\sigma)\simeq
\Bbb Z_3\ltimes
\Bbb Z_2^3.$  Let
$s=(123).$  Let
$\kappa$ be a character of $R(\pi),$ and let $R_\kappa$ be the stabilizer of $\kappa$ in $R(\sigma).$  By Proposition 25 of
\cite{\bf Se} we know that any irreducible representation of $R(\sigma)$ is given by
$\rho=\rho_{\kappa,\lambda}=\Ind_{R_\kappa}^{R(\sigma)}(\kappa\otimes\lambda),$ where $\kappa$ is extended to $R_\kappa$
trivially, and $\lambda$ is an irreducible representation of $R_{\kappa}\cap \Bbb Z_3.$  
 Note
that if
$w=sc,$ then $w$ acts transitively on $R(\pi).$   Thus $R_\kappa\neq R(\pi)$ if and only if $\kappa(C_i)=\kappa(C_j)$
for $i\neq j.$  This implies either $\kappa=1,$ or $\kappa(C_i)=-1$ for all $i,$   and we denote this character  by $sgn.$ 
So, if $\kappa\neq1,sgn,$ then $R_\kappa=R(\pi),$ and $\rho=\rho_\kappa=\Ind_{R(\pi)}^{R(\sigma)}\kappa.$
Then, by the induced character formula, the character $\theta_{\kappa}$ of $\rho_\kappa$ on an element $w=sc$ is
$$\theta_\kappa(sc)= \frac18\sum_{\matrix x\in R(\sigma)\\x^{-1}(sc)x\in R(\pi)\endmatrix}\kappa(x^{-1}scx)=0,$$
as the sum has no terms since $R(\pi)$ is normal.  Thus,  all the components $\pi(\rho_{\kappa})$ with $\kappa\neq
1,sgn$ are non-elliptic. There are two such representations, as there are two orbits of such $\kappa$ under $\Bbb
Z_3.$ Each of these components appears in $i_{G,M}(\sigma)$ with multiplicity three.  On the other
hand, the six components components
$\pi(\rho_{1,\lambda}),\,\pi(\rho_{sgn,\lambda}),$ are elliptic, and these appear with multiplicity one.  We conclude
that there are representations $\sigma$ for which
$i_{G,M}(\sigma)$ has some elliptic components, but not all the components are elliptic.  It is clear that this
phenomenon will generalize to the case where
$r$ is prime.

Note that, when $P=B$ is the Borel subgroup, these are precisely the examples discussed in \cite{\bf K2}.  The phenomenon noted above
was not mentioned there merely because the results of \cite{\bf A2} were not yet available.

\end